\documentclass{article}
\usepackage{amsmath}
\usepackage{amsfonts}
\usepackage{amssymb}
\usepackage{graphicx}

\pdfoutput=1
\newtheorem{theorem}{Theorem}

\newtheorem{corollary}[theorem]{Corollary}

\newtheorem{lemma}[theorem]{Lemma}

\newenvironment{proof}[1][Proof]{\noindent\textbf{#1.} }{\ \rule{0.5em}{0.5em}}
\input{tcilatex}

\begin{document}

\title{On the method of likelihood-induced priors\thanks{%
Presented at MaxEnt 2018, 38th International Workshop on Bayesian Inference
and Maximum Entropy Methods in Science and Engineering (July 02-06, 2018,
The Alan Turing institute , London, UK).}}
\author{Ali Ghaderi\thanks{%
E-mail: Ali.Ghaderi@usn.no} \\
%EndAName
University of South-Eastern Norway\\
Kj\o lnes Ring 56, NO-3901 Porsgrunn, Norway}
\maketitle

\begin{abstract}
We demonstrate that the functional form of the likelihood contains a
sufficient amount of information for constructing a prior for the unknown
parameters. We develop a four-step algorithm by invoking the information
entropy as the measure of uncertainty and show how the information gained
from coarse-graining and resolving power of the likelihood can be used to
construct the likelihood-induced priors. As a consequence, we show that if
the data model density belongs to the exponential family, the
likelihood-induced prior is the conjugate prior to the corresponding
likelihood. 
\end{abstract}

\section{Introduction}

We argue that the functional form of the likelihood function is informative
and in the absence of other types of information, it induces a prior on the
unknown parameters of interest. Regardless of the data, the functional form
of the likelihood depends on the design of experiment, the measurement
methods and the model which is evaluated. In the following we demonstrate
how this type of information can be used to construct priors for the
parameters of interest.

\section{Statement of the problem}

In the Bayesian approach, the possible values of the unknown, say $\mathbf{%
\theta }$ are described by the posterior distribution $p\left( \pmb{\theta }%
\left\vert \pmb{d}\right. \right) $ through the Bayes' theorem%
\begin{equation}
p\left( \pmb{\theta }\left\vert \pmb{d}\right. \right) =\frac{p\left( \pmb{d}%
\left\vert \pmb{\theta }\right. \right) p(\pmb{\theta })}{p\left( \pmb{d}%
\right) }
\end{equation}%
where $p\left( \pmb{d}\left\vert \pmb{\theta }\right. \right) $ is the
likelihood, $p(\pmb{\theta })$ is the prior and $p\left( \pmb{d}\right) $ is
the evidence, also known as the marginal likelihood. The likelihood is often
well-determined by the process model and the data-generating process \cite%
{Gregory2005}. The prior encodes ones belief about $\pmb{\theta }$ before
seeing the data. One of the challenges is to construct a prior that reflects
ones state of knowledge before seeing the data. In many applications there
are often little or no prior information available about $\pmb{\theta .}$
Sometimes it is even difficult to interpret $\pmb{\theta ,}$ let alone
describing the degree of belief. In such cases, there are several existing
methods for constructing noninformative priors.

\begin{description}
\item[Uniform prior] Lack of information is modelled by \emph{bounded}
uniform distribution \cite{Laplace1814}.

\item[Jeffreys prior] The prior is proportional to square root of the fisher
information derived from the likelihood function \cite{Jeffreys1998}.

\item[Reference prior] The idea is to obtain a prior that maximizes the
expected gain of information provided by the data \cite{berger2009}.

\item[Maximal data information prior] The approach is based on the
information conservation principle \cite{Zellner1971}.
\end{description}

Although each one of these methods have their domain of usefulness, they
implicitly ignore certain aspects of the information in a crucial way. The
functional form of the likelihood contains some important information. In
this sense, we are not ignorant about $\pmb{\theta .}$ In the following, we
demonstrate how the knowledge about the functional form of the likelihood
imposes some sort of constraint on our belief, which in turn induces a prior
on $\pmb{\theta }$.

\section{Measure of uncertainty}

Uncertainty is a direct result of having multitude of choices. The
assignment of probabilities to each choice reflects our belief. The
certainty is only achieved in the limit when the number of choices is
reduced to one. Therefore, in statistical inference, one aspires to achieve
greater certainty by aggregation of choices. The purpose of the aggregation
of any quantity is to create a less uncertain description without loss of
essential information. In general, if $\pmb{\theta }$ is the unknown, the
information preserving aggregation is applied to possible values /functions
of $\pmb{\theta }$. Following this path of thought, two related questions
present themselves,

\begin{enumerate}
\item What kind of information is preserved?

\item What does such information tell us about the probabilities of the
events of interest?
\end{enumerate}

In the following, we shall call the application of any information
preserving aggregation as \emph{coarse-graining}. In the probabilistic
inference, the coarse-graining is conducted by applying the expectation
operator with respect to a probability distribution on the quantity of
interest. Moreover, since we are only interested on the coarse grainable
quantities, such quantities will be referred to as \emph{conservation laws}.
In general, aggregation is destructive. If the aggregation is to be
conservative then the distribution of the quantity of interest has to
contain a sufficient amount of redundancy. More generally, regular
structures contain a lot of redundancies whilst the random structures
contain none. In the latter case, the aggregation would result in loss of
information.

Since coarse-graining reduces the uncertainty, any study of such processes
would depend on our ability to quantify uncertainty. That is, since
uncertainty about a quantity is modelled by a probability distribution then
we need a way to measure the amount of uncertainty described by a given
probability distribution. This measure has already been defined in the
closely related field of the information theory. Indeed, the fact that in
statistical inference it is often possible to make inferences on the
coarse-grained information bears the resemblance to the goal of information
theory, which is to compress information such that it can be recovered
exactly or approximately. In the information theory, Claude Shannon
described three requirements which any measure of uncertainty should
fulfill. However, in the present context, these requirements need a slight
modification in order to meet our needs. Accordingly, these requirements can
be stated as follows

\begin{description}
\item[Continuity:] The measure of uncertainty is a continuous function of
probabilities. This means that small changes in the value of probabilities
should only change the measure by a small amount.

\item[Maximum] Maximum uncertainty is achieved when for all $i$ we have $%
P_{i}=M_{i},$ where $M_{i}$'s are determined by the resolving power of the
available information.

\item[Conservation of uncertainty:] Uncertainty in the fine
grained-description is equal to the sum of uncertainties in the
coarse-grained description and the amount concealed in each grain.
\end{description}

The maximum uncertainty is equivalent to lack of redundancy or compressible
structure. That is, at the state of maximum uncertainty, aggregation of
events would result in loss of information. At this state, the only thing
one can do is to assign weights to events based on the resolving power of
the measurements. This is a reasonable requirement in applications where
measurements are part of inference. In the context of the transmission of
messages, Shannon assumed that $M$ is the uniform distribution (see \cite%
{lesne2014} for further insight on Shannon's entropy and related subjects
relevant to the present topic). The third requirement, the conservation of
uncertainty, follows from the product and sum rules which impose constraints
on how the probabilities of the events at the fine-grained level scale under
coarse-graining. The functional that fulfils these three conditions is known
as \emph{entropy} and in the discrete case is given as%
\begin{equation}
H_{p}=-\sum_{i=1}^{n}P(x_{i})\ln \frac{P(x_{i})}{M(x_{i})}.
\label{entropy-d}
\end{equation}%
This can be extended in the limit to the continuous case%
\begin{equation}
H_{p}=-\int_{\Omega }p(x)\ln \frac{p(x)}{m(x)}dx  \label{entropy-c}
\end{equation}%
where $\Omega $ denotes the support of $p$ (see also \cite{Shore1980} for
other axiomatic derivation). In order for this expression to be
well-defined, $m$ has to dominate $p,$ that is, for almost all $x$ if $%
m(x)=0 $ then $p(x)=0.$ In the discrete case, the positive function $M$ is
the \emph{counting measure} and in the continuous case $m$ is the \emph{%
Lebesgue measure}. These measures can also be thought as mechanisms for
assigning equal mass to regions of equal volume, which in the continuous
case, guaranties that $H$ is invariant with respect to change of variables.

If $m$ is normalized then it follows from Jensen's inequality that $H\leq 0,$
where the equality is achieved if and only if $p=m.$ Thus the uncertainty
described by $p(x)$ is ranked relative to $m(x).$ In the following, we will
only consider the continuous case ( \ref{entropy-c}) and it is assumed that $%
m$ is a probability density and all the integrals are proper and finite. The
improper integrals are considered as the limit of proper integrals and are
intended as useful approximations. In such cases, there is no guaranty that $%
m$ is normalizable and hence the upper limit of $H$ might be larger than
zero or even without limit.

\section{Conservation laws}

The maximum likelihood method plays an important role in frequentist
inference. The uncertainty reduction is conducted by choosing $\theta $ to
be the global maximum of the likelihood. Since the logarithm is a
monotonically increasing function, the maximum of the likelihood function is
the same as the maximum of the corresponding average log-likelihood. For
independent and identically distributed observations drawn from the data
model density $p\left( x\left\vert \theta \right. \right) ,$ the likelihood
and the average log-likelihood are defined as%
\begin{equation}
L(\theta ;x_{1:n})=\prod\limits_{i=1}^{n}p\left( \left. x_{i}\right\vert
\theta \right)
\end{equation}%
and%
\begin{equation}
l(\theta ;x_{1:n})=\frac{1}{n}\ln \left\{ L(\theta ;x_{1:n})\right\} =\frac{1%
}{n}\sum\limits_{i=1}^{n}\ln \left\{ p\left( \left. x_{i}\right\vert \theta
\right) \right\} ,
\end{equation}%
respectively. Then the maximum likelihood estimate (MLE) is%
\begin{equation}
\widehat{\theta }=\arg \underset{\theta }{\max }L(\theta ;x_{1:n})=\arg 
\underset{\theta }{\max }l(\theta ;x_{1:n}).
\end{equation}

The use of MLE is justified in the limit when the number of observations
grows to infinity. However, in reality we have only a finite number of
observations. Therefore, the true maximum of the log-likelihood is not known
and hence, in general, we can only talk about plausible candidates for $%
\theta .$ In the Bayesian context, the uncertainty about the average
log-likelihood, $l(\theta ;x_{1:n}),$ can be modelled by a probability
distribution. For a given finite set of observations, $l(\theta ;x_{1:n})$
depends only on $\theta $. Since we do not know the true value of $\theta ,$
by appropriate coarse-graining with respect to $\theta ,$ we can say
something about the centre of mass of the distribution for $l(\theta
;x_{1:n}).$ Often this coarse-graining reveals conservation laws, which as
we shall see shortly, can be used to identify the family of distributions
describing the uncertainty in $\theta .$

Let's demonstrate this by an example in which the observations are generated
from an exponential distribution with the parameter $\theta .$ Assume that
the observations consists of $n$ independent data points, $x_{1:n}.$ Then
the likelihood is given by%
\begin{equation}
L(\theta ;x_{1:n})=\prod\limits_{i=1}^{n}p\left( x_{i}\left\vert \theta
\right. \right) =\frac{1}{\theta ^{n}}\exp \left( -\frac{1}{\theta }%
\sum_{i=1}^{n}x_{i}\right)
\end{equation}%
and the average log-likelihood is%
\begin{equation}
l(\theta ;x_{1:n})=-\ln \theta -\theta ^{-1}\bar{x},
\end{equation}%
where $\bar{x}$ is the arithmetic average of the observations. Clearly,
uncertainty in $\theta $ results in uncertainty in $l(\theta ;x_{1:n})$.
Coarse-graining of $l(\theta ;x_{1:n})$ with respect to $\theta $ results in
a single number, which is the centre of mass for the distribution assigned
to $l(\theta ;x_{1:n}).$ The coarse-graining with respect to $\theta $ is
conducted by taking the expectation with respect to $p\left( \theta
|x_{1:n}\right) ,$ i.e.%
\begin{equation}
E_{\theta |\mathbf{x}}\left[ l(\theta ;x_{1:n})\right] =-\left\langle \ln
\theta \right\rangle -\left\langle \theta ^{-1}\right\rangle \bar{x}
\end{equation}%
where the operator $\left\langle \cdot \right\rangle $ denotes the
expectation taken with respect to $p\left( \theta |x_{1:n}\right) .$ This
closely follows the Bayesian philosophy that if one does not know the true
value of $\theta $ then one should average over all its possible values. In
the present example, it is evident that all the information from the
distribution of $\theta $ relevant for determining $\left\langle l(\theta
;x_{1:n})\right\rangle $ are summarized by the numbers $\left\langle \ln
\theta \right\rangle $ and $\left\langle \theta ^{-1}\right\rangle .$ This
implies that the information about $E\left( l\right) $ induces a class of
distributions for $\theta $ that conserve the expected values of $%
f_{1}(\theta )=\ln \theta $ and $f_{2}(\theta )=\theta ^{-1}$ consistent
with $E\left( l\right) .$ In general, this class is very large. However, as
we shall see, in the absence of any other information, these likelihood
induced conservation laws can be used to identify the parametric family that 
$p\left( \theta |x_{1:n}\right) $ belongs to.

\section{The resolving power}

In practice, the resolving power of any measurement system is finite. That
is, if $\theta _{1}$ and $\theta _{2}$ are too close or similar, one may not
expect to detain decisive support for $\theta _{1}$ against $\theta _{2}$
from the data. Often the resolving power differs from one region of the
parameter space to another. Taking this into account, in inference, we
should give more weight to the regions for which we can potentially detain
strong support from the data. As we will demonstrate shortly, this weighting
procedure will result in a density, which we shall define as the state of
maximum uncertainty.

Although, a priori we do not know the data, nevertheless, the data model
density $p\left( x|\theta \right) $ contains some information about the
resolving power of the data-generating process. One possible approach is to
consider the sensitivity of $p\left( x|\theta \right) $ to the changes in $%
\theta .$ To this end, the score function (the derivative of the $p\left(
x|\theta \right) $ normalized by its value) is a good indication of the
sensitivity. It is given by%
\begin{equation}
S\left( \theta ,X\right) =\frac{1}{p\left( x|\theta \right) }\frac{\partial
p\left( x|\theta \right) }{\partial \theta }=\frac{\partial \ln p\left(
x|\theta \right) }{\partial \theta }.
\end{equation}%
We are interested in the score as a function of $\theta .$ For given $X,$
large absolute values of the score indicate high sensitivity and hence high
resolving power. Since a priori the data is not known, the true value of the
score function is uncertain. Nonetheless, we can consider the mean and the
variance of the score as an indication of the resolving power. It can be
shown that under some regularity condition, the mean of score function is%
\begin{equation}
E\left( S|\theta \right) =\int_{\Omega }\frac{\partial \ln p\left( x|\theta
\right) }{\partial \theta }p\left( x|\theta \right) dx=0
\end{equation}%
and its variance is%
\begin{equation}
\mathcal{I}\left( \theta \right) =\int_{\Omega }\left( \frac{\partial \ln
p\left( x|\theta \right) }{\partial \theta }\right) ^{2}p\left( x|\theta
\right) dx\geq 0
\end{equation}%
where $\mathcal{I}\left( \theta \right) $ is also known as the \emph{Fisher
information}. The larger the Fisher information is, the greater is the
chance of observing score values at larger distance from zero. Thus, for a
given $\theta ,$ the large value of the Fisher information indicate high
chance of having high resolving power in the neighbourhood of that specific $%
\theta .$ This observation suggests that the Fisher information can be used
for weighting $\theta $ according to its probable degree of resolution.
However, at least two problems present themselves:

\begin{enumerate}
\item In general, due to multidimensionality of the parameter $\theta ,$ the
Fisher information is a matrix and not a scalar.

\item The dimension of the Fisher information is $\left[ \theta ^{-2}\right]
.$ It needs to be $\left[ \theta ^{-1}\right] $ in order for $m\left( \theta
\right) d\theta $ to be dimensionless.
\end{enumerate}

The first problem implies that our resolution might depend on the direction
we move. Nevertheless, the volume of the $n$-parallelotope spanned by the
column vectors of the Fisher information matrix can be used as an indication
of the resolving power. This volume is found by taking the determinant of $%
\mathcal{I}\left( \theta \right) .$ The second problem can be addressed by
taking the square root of the volume. Consequently, the density that
describes the maximum uncertainty is%
\begin{equation}
m(\theta )d\theta \propto \sqrt{\det \mathcal{I}\left( \theta \right) }%
d\theta .  \label{Jeffreys prior}
\end{equation}%
In literature, this density is known as the \emph{Jeffreys prior}. This
result is a direct consequence of our way of assigning prior weights to
regions of parameter space with respect to probable resolving power of the
likelihood.

\section{Induced priors}

In the previous sections, we demonstrated that by coarse-graining of average
log-likelihood, one can identify the essential information from the
distribution of $\theta $ relevant to the centre of mass of the average
log-likelihood distribution. This type of information, if attainable, come
as a set of conservation laws. We have also argued that Jeffreys prior, eq. (%
\ref{Jeffreys prior}), can serve as the density with the maximum
uncertainty. In the absence of any other information except the functional
form of the likelihood, the question about the prior on $\theta $ can be
formulated as follows: If we had previously seen $n$ observations, what can
we say about the value of $\theta $ before seeing the new observations? In
order to be able to answer this question we need the following Lemma and its
corollary.

\begin{lemma}
\label{Lem-MLEnt}Let $q(\theta |x)$ dominate $p(\theta |x).$ Then%
\begin{equation*}
-\int_{\Omega }p(\theta |x)\ln \left( \frac{p(\theta |x)}{q(\theta |x)}%
\right) d\theta \leq 0
\end{equation*}

\begin{proof}
The result follows from Jensen's inequality.
\end{proof}
\end{lemma}

\begin{corollary}
\label{Cor-MLEnt}Let $q(\theta |x)$ dominate $p(\theta |x)$ with respect to
the common measure $m(\theta ).$ Then%
\begin{equation}
-\int_{\Omega }p(\theta |x)\ln \left( \frac{q(\theta |x)}{m(\theta )}\right)
d\theta \geq -\int_{\Omega }p(\theta |x)\ln \left( \frac{p(\theta |x)}{%
m(\theta )}\right) d\theta =H_{p}(x;m)  \label{GibbsInequality}
\end{equation}

\begin{proof}
The statement follows from lemma \ref{Lem-MLEnt}.
\end{proof}
\end{corollary}

The statement of the above corollary is the same as the Gibbs' inequality
with respect to the common measure $m(\theta ).$ Further, it is assumed that 
$m(\theta )$ is a density which dominates both $q(\theta |x)$ and $p(\theta
|x).$

Now, let $f_{k}(\theta )$ be $r$ different functions of $\theta $ such that 
\begin{equation}
\int_{\Omega }f_{k}(\theta )p(\theta |x)d\theta =\int_{\Omega }f_{k}(\theta
)q(\theta |x)d\theta =F_{k},\text{ for all }k=1,\ldots ,r  \label{ExpctCons}
\end{equation}%
In general, there are many densities which satisfy the above constraints.
However, as it will become clear further below, we are interested on the
densities with maximum entropy. The following is an extension of the
argument given by Jaynes for the discrete case \cite[p.357]{Jaynes2003}. Let 
\begin{equation}
q(\theta |x)=\frac{m(\theta )}{Z\left( \lambda _{1},\ldots ,\lambda
_{r}\right) }\exp \left( -\sum_{k=1}^{r}\lambda _{k}f_{k}\left( \theta
\right) \right)  \label{ExpoDistQ}
\end{equation}%
where $Z\left( \lambda _{1},\ldots ,\lambda _{r}\right) $ is the
normalization constant and $\lambda _{k}$ are functions of $x.$ It can be
shown, by the method of Lagrange multipliers, that $q$ satisfies the
constraints (\ref{ExpctCons}). After substituting $q$ into eq. (\ref%
{GibbsInequality}) and taking into account the constraints (\ref{ExpctCons}%
), we get%
\begin{equation}
H_{p}(x;m)=-\int_{\Omega }p(\theta |x)\ln \left( \frac{p(\theta |x)}{%
m(\theta )}\right) d\theta \leq \ln Z+\sum_{k=1}^{r}\lambda _{k}F_{k}
\end{equation}%
This relation holds for all $p$ which are dominated by $q.$ In particular,
this inequality holds for the class of all the densities which satisfy the
constraints (\ref{ExpctCons}) and are dominated by $q$. The equality is only
achieved if $p=q.$ In effect, it can be seen that every density with respect
to the measure $m$ has lower entropy than $q$. Indeed, notice that all the
zeros of $q$ coincide with the zeros of $m.$ Since $m$ dominates all the
densities in the class of interest, then we can conclude that $q$ also
dominates every density which is dominated by $m$ and hence $q$ dominates
every density in the class of interest.

At this point one might raise the following question that why should one
pick the density with the maximum entropy. Recall, that the entropy can also
be interpreted as a measure of redundancy. Any other family of densities
satisfying the constraints (\ref{ExpctCons}) have higher redundancy, and
hence possibly, larger number of conservation laws, some of which are not
induced by the likelihood function.

In the above the numerical values of the observations were irrelevant.
Therefore, we can consider them as the pseudo-observations. We are now ready
to lay out the algorithm for constructing the likelihood-induced prior. The
algorithm is as follows

\begin{enumerate}
\item Determine the average Log-likelihood using the pseudo-observations $%
x_{1:n}$.

\item Conduct coarse-graining with respect to $p\left( \theta
|x_{1:n}\right) $ and identify the conservation laws.

\item Determine the Jeffreys' prior from the data model $p\left( \theta
|x\right) $.

\item Construct the maximum entropy density given by eq. (\ref{ExpoDistQ})
\end{enumerate}

Let us demonstrate this algorithm for two cases.

\subsection{Exponential density}

Previously, we described the conservation laws for the data model density
being exponential with respect to $x.$ Below we list the results after each
step of the algorithm without details.

\begin{enumerate}
\item Average log-likelihood: $l=n^{-1}\ln L=-\ln \theta -\theta ^{-1}%
\overline{\pmb{x}}.$

\item Coarse-graining: $E_{\theta }\left( l\right) =-\left\langle \ln \theta
\right\rangle -\left\langle \theta ^{-1}\right\rangle \overline{\pmb{x}}%
\Rightarrow $ the conservation laws are $f_{1}\left( \theta \right) =\ln
\theta $ and $f_{2}\left( \theta \right) =\theta ^{-1}.$

\item Jefreys' prior: $m\left( \theta \right) \propto \theta ^{-1}$
\end{enumerate}

Applying the step $4$ of the algorithm results in the parametric family of
inverse-gamma distributions%
\begin{equation}
q\left( \theta \right) =\frac{1}{\lambda _{2}\Gamma \left( \lambda
_{1}\right) }\left( \frac{\theta }{\lambda _{2}}\right) ^{-\lambda
_{1}-1}\exp \left( -\frac{\lambda _{2}}{\theta }\right) ,\text{ where }%
\lambda _{1},\lambda _{2}>0
\end{equation}

\subsection{Exponential family}

A density belongs to the exponential family if it can be expressed in the
following way%
\begin{equation}
p\left( \pmb{x}\left\vert \pmb{\eta }\right. \right) =h\left( \pmb{x}\right)
\exp \left\{ \pmb{\eta }^{T}\cdot T\left( \pmb{x}\right) \right\} g\left( %
\pmb{\eta }\right)
\end{equation}%
where $\pmb{\eta }=$ $\pmb{\eta }_{1:s}\left( \pmb{\theta }\right) $ is the
natural parameters, $\pmb{\eta }^{T}\cdot T\left( \pmb{x}\right)
=\sum_{i=1}^{s}\eta _{i}T_{i}(\pmb{x})$ and the support of the density does
not depend on the choice of $\pmb{\theta }$. If the data model density
belongs to the exponential family, then the likelihood of the $n$ iid
observations $\pmb{x}_{1:n}$ is%
\begin{equation}
L\left( \pmb{\theta };\pmb{x}_{1:n}\right) =\prod\limits_{k=1}^{n}p\left(
x_{k}\left\vert \pmb{\eta }\right. \right) =\left(
\prod\limits_{k=1}^{n}h\left( x_{k}\right) \right) g\left( \pmb{\eta }%
\right) ^{n}\exp \left\{ \pmb{\eta }^{T}\cdot T\left( \pmb{x}\right)
\right\} .
\end{equation}%
Below we list the results of each step of the algorithm without details.

\begin{enumerate}
\item Average log-likelihood: $l=n^{-1}\ln L=\overline{\ln h\left( \pmb{x}%
\right) }+\ln g\left( \pmb{\eta }\right) +\pmb{\eta }^{T}\cdot \overline{%
T\left( \pmb{x}\right) }$.

\item Coarse-graining: $E_{\pmb{\theta }}\left( l\right) =\overline{\ln
h\left( \pmb{x}\right) }+\left\langle \ln g\left( \pmb{\eta }\right)
\right\rangle +\left\langle \pmb{\eta }\right\rangle ^{T}\cdot \overline{%
T\left( \pmb{x}\right) }\Rightarrow $ the conservation laws are $f_{1}\left( %
\pmb{\theta }\right) =\ln g\left( \pmb{\eta }\right) $ and $f_{2}\left( %
\pmb{\theta }\right) =\pmb{\eta }.$

\item Jeffreys prior: $m\left( \pmb{\theta }\right) \propto \sqrt{\det 
\mathcal{I}\left( \pmb{\theta }\right) }$
\end{enumerate}

Applying the final step of the algorithm results in%
\begin{equation}
q\left( \pmb{\theta }\right) \propto \sqrt{\det \mathcal{I}\left( %
\pmb{\theta }\right) }g\left( \pmb{\eta }\right) ^{\gamma }\exp \left( %
\pmb{\eta }^{T}\cdot \pmb{\lambda }\right) .
\end{equation}%
In literature, the distribution $q$ is known as the \emph{conjugate prior}
for the likelihood $L.$ The conjugate priors play an important role in
Bayesian statistics. They are often used because they result in closed form
expression for posterior without the need for calculation of the elusive
normalization constant. However, the result of this section further
elucidates the role of the conjugate priors as not only algebraically
convenient constructs but as manifestation of the likelihood-induced
information.

\section{Concluding remarks}

We have demonstrated that, in general, the functional form of the likelihood
contains enough information for constructing a prior for the unknown
parameters. The key idea was the coarse-graining of the information in order
to reduce uncertainty and using entropy as the measure of uncertainty. The
identification of the relevant conservation laws along with the measure for
resolving power resulted in a four-step algorithm for constructing the
likelihood-induced priors. We have also demonstrated that in case the data
model density belongs to the exponential family, the likelihood-induced
prior is the conjugate prior for the corresponding likelihood. Furthermore,
this algorithm can readily be applied to other parametric classes of
densities. We shall come back to this issue in the future.

\bibliographystyle{alpha}
\bibliography{CitiationDB}

\end{document}